\RequirePackage[online]{threeparttable}

\newif\ifSEA
\SEAfalse

\ifSEA
	\documentclass[a4paper,UKenglish,cleveref, autoref, thm-restate]{lipics-v2021}
\else
	\documentclass[notitlepage,11pt]{article}

	\usepackage[left=1.25in,right=1in,top=1in,bottom=1in]{geometry}
	\usepackage{amsmath,amsthm,amssymb,amsfonts}
	\usepackage{thmtools, thm-restate}
	\usepackage[title,header]{appendix}
	\usepackage{hyperref}
	\hypersetup{colorlinks=true,citecolor=blue,linkcolor=red,urlcolor=blue,breaklinks=true}	
	
	\declaretheorem{theorem}
	\declaretheorem{lemma, proposition, corollary}[sibling=theorem]
	\declaretheorem{definition, assumption, conjecture, fact, question}[style=definition]
	
	\usepackage[maxbibnames=99,maxcitenames=2,backend=biber,sorting=nyt,sortcites,bibencoding=utf8,natbib=true]{biblatex}
	\AtEveryBibitem{
		\ifentrytype{article}{\clearfield{url}\clearfield{language}\clearfield{publisher}}{}
		\ifentrytype{book}{\clearfield{doi}}{\clearfield{publisher}}
		\ifentrytype{incollection}{\clearfield{doi}}{}
		\ifentrytype{inproceedings}{\clearfield{doi}}{}
		\iffieldundef{volume}{}{\clearfield{doi}\clearfield{eprint}}
		\ifentrytype{unpublished}{\iflistundef{publisher}{}{\clearfield{url}}}{\clearfield{url}}
		\clearlist{language}
		\clearfield{isbn}
		\clearfield{issn}
		\clearfield{month}
	}
	
	\usepackage[nameinlink]{cleveref}
	\crefname{equation}{}{}%{eq.}{eqs.}
	\crefname{section}{\textsection}{\textsection}
	\AtBeginEnvironment{appendices}{\crefalias{section}{appendix}\crefalias{subsection}{appendix}}
	\Crefname{secinapp}{Appendix}{Appendices}
	
	\renewcommand{\ref}[1]{\cref{#1}}
\fi

\usepackage{bm,booktabs}
\usepackage{adjustbox}

\newlength{\LPlength}
\newcommand{\LPblocktag}[2]{\settowidth{\LPlength}{#1}%
                            \parbox{\LPlength}{\begin{align}\tag{#1}#2\end{align}}%
                            \hspace*{\fill}}

\title{Multilinear formulations for computing a Nash equilibrium of multi-player games}

\ifSEA
	\titlerunning{Nash equilibria through multilinear optimization}
	\author{Miriam Fischer}{Department of Computing, Imperial College London, United Kingdom \and \url{https://www.doc.ic.ac.uk/~mif21/}}{m.fischer21@imperial.ac.uk}{}{Supported by a PhD scholarship from DeepMind}
	\author{Akshay Gupte\footnote{Corresponding author}}{School of Mathematics, The University of Edinburgh, United Kingdom \and \url{https://www.maths.ed.ac.uk/school-of-mathematics/people/a-z?person=820}}{akshay.gupte@ed.ac.uk}{https://orcid.org/
0000-0002-7839-165X}{}
	\authorrunning{M. Fischer and A. Gupte}
	\Copyright{Miriam Fischer and Akshay Gupte} %TODO mandatory, please use full first names. LIPIcs license is "CC-BY";  http://creativecommons.org/licenses/by/3.0/


 %TODO mandatory; please add comma-separated list of keywords

\category{} %optional, e.g. invited paper

\acknowledgements{This research was initiated as part of the MSc dissertation of the first author in the School of Mathematics at the University of Edinburgh}%optional

\supplementdetails[linktext={GitHub code}]{Instances of games and source code for optimization models}{https://github.com/economicsandcomputing/MultilinearNashEquilibria}
\else
	\author{Miriam Fischer \thanks{Department of Computing, Imperial College London, UK. \texttt{m.fischer21@imperial.ac.uk}. Research was initiated as part of MSc dissertation in the School of Mathematics at the University of Edinburgh} \textsuperscript{} \thanks{Supported by a PhD scholarship from DeepMind} \and Akshay Gupte \thanks{School of Mathematics, The University of Edinburgh, UK. \texttt{akshay.gupte@ed.ac.uk}}}
	\date{Accepted: April 2023}
\fi

\begin{document}
\maketitle

\begin{abstract}
We present multilinear and mixed-integer multilinear programs to find a Nash equilibrium in multi-player noncooperative games. We compare the formulations to common algorithms in Gambit, and conclude that a multilinear feasibility program finds a Nash equilibrium faster than any of the methods we compare it to, including the quantal response equilibrium method, which is recommended for large games. Hence, the multilinear feasibility program is an alternative method to find a Nash equilibrium in multi-player games, and outperforms many common algorithms. The mixed-integer formulations are generalisations of known mixed-integer programs for two-player games, however unlike two-player games, these mixed-integer programs do not give better performance than existing algorithms.
\ifSEA
\else
\fi
\end{abstract}

\section{Introduction}
A noncooperative game has $n$ players, where $n\ge 2$ is finite, with each player having finitely many pure strategies which they do not discuss or reveal to each other. A mixed strategy for a player is a probability distribution over the player's pure strategies. Each player has a known payoff function which maps any combination of pure strategies of all the $n$ players to a real number. Mixed strategies of all the players form a tuple whose payoff is calculated by taking expectation over the probability distributions. In his seminal work \cite{nash1951}, Nash showed that every such game has a tuple of mixed strategies that is an equilibrium in the sense that no player increases their payoff if they were to change their mixed strategy while the others keep theirs fixed. Although existence of Nash equilibrium is guaranteed, uniqueness does not always hold, and there are also characterisations of when there exists an equilibrium formed solely by pure strategies. 

This paper deals with the question of algorithmic and numerical computation of Nash equilibria. From a complexity perspective, computing an equilibrium was only somewhat recently formally settled to being PPAD-complete \cite{chen2009settling,daskalakis2009} even for two-player games. There is a lot of literature for two-player bimatrix games, and the most well-known and established exact method to compute an equilibrium is the Lemke-Howson algorithm \cite{Howson1964}. This gives a very good computational performance on many instances in practice, although its worst-case performance can take exponentially many pivoting steps \cite{Stengel2006}. 

However, for multi-player games, there do not seem to be commonly established approach for computing the equilibrium. Although there is a generalisation of the Lemke-Howson method to $n$-person games \cite{Rosenmueller1971,wilson1971}, popular algorithmic approaches include a global Newton method \cite{Govindan2003}, an iterated polymatrix approximation approach \cite{Govindan2004}, a simplicial subdivision method \cite{Laan1987}, a simple search algorithm aiming to find an equilibrium with small support size \cite{Porter2008}, and a quantal response equilibrium method which gives an approximation to a Nash equilibrium \cite{Turocy2005}. Many of these methods are implemented in the game-theoretic library \texttt{Gambit} \cite{McKelvey2014}. Experiments comparing different methods have been undertaken \cite{Sandholm2005, Porter2008}, however it is rather unclear which of the methods is best for multi-player games. For example, the global Newton method gives solid performance for small games, however does not to scale well to larger games \cite{computation2021,TurocyAlgorithm}. Support enumeration algorithms are fast for games with pure equilibria but will be much slower for a game that only has equilibria of medium to large support size. There are also approximation algorithms, which tend to approximate a Nash equilibrium for large games \cite{Turocy2005,computation2021,gamesapproximation2016}.

We adopt the optimization approach, and propose different optimization formulations for computing a Nash equilibrium for $n$-person games for $n \geq 2$. Particularly, we present a multilinear polynomial continuous feasibility program of degree $n$ (= number of players), which is a generalisation of the bilinear optimization problem for 2 players \cite{Mangasarian1964}. Further, we extend the two-player mixed-integer formulations of \cite{Sandholm2005} to multi-player games, and give a large variety of mixed-integer formulations to find a Nash equilibrium in multi-player games. All our formulations find a Nash equilibrium of a $n \geq 2$ player game. We compare our programs to \texttt{gambit-gnm} (global Newton), \texttt{gambit-simpdiv} (simplicial subdivision), \texttt{gambit-logit} (quantal response equilibrium) algorithms in \texttt{Gambit}, focusing on random games and covariant games with negative covariance. We find that the mixed-integer formulations do not give better performance than existing algorithms, and our analysis of those is aimed to get an understanding of which mixed-integer formulations are most suited for finding a Nash equilibrium. We find that our multilinear continuous feasibility program is faster than all the methods in Gambit we compare it to, including the \texttt{gambit-logit} method, which is so far recommended for large games. Thus, we provide an alternative approach to computing Nash equilibrium in multi-player games. 

The next section presents our continuous and mixed-integer multilinear optimization formulations. For each of them, their correctness, i.e., the fact that their optimal/feasible solutions correspond to Nash equilibria of the game, is proved in the Appendix. 

\section{Formulations} \label{sec:formulation}

The multi-player multilinear formulation is an extension of a bilinear formulation for bimatrix games \cite{Mangasarian1964}. To motivate the multilinear formulation, we shortly recall the bilinear program that is equivalent to finding a Nash equilibrium in a bimatrix game. To do so, we introduce some notation. Let $A,B \in \mathbb{R}^{m \times n}$ be the payoff matrix of player 1 and player 2, with $m$ pure strategies of player 1 and $n$ pure strategies of player 2. Let $\bm{x} \in \mathbb{R}^m$ with $\bm{x} \geq 0$ and $\sum_{i=1}^m x_i = 1$ be a (possibly mixed) strategy of player 1, with $x_s$ being the probability placed on pure strategy $s$. Let $\bm{y} \in \mathbb{R}^n$ with $\bm{y} \geq 0$ and $\sum_{j=1}^n y_j = 1$ be a (possibly mixed) strategy of player 2. Let $\bm{1}_n$ and $\bm{1}_m$ denote vectors of all ones of dimension $n$ and $m$. Any globally \textit{optimal} solution $(x,y,p,q)$ to the bilinear optimization problem in \ref{fig:bimatrixform} is equivalent to a Nash equilibrium in a bimatrix game. 

\LPblocktag{BLP}{\label{fig:bimatrixform}}%
\begin{minipage}{12cm}
\begin{subequations}
\begin{align}
\max_{x,y,p,q} &\quad x^\top A y + x^\top B y - p - q \\
\text{s.t.} &\quad Ay \le p\bm{1}_m, \quad B^{\top}x \le  q\bm{1}_n\\
&\quad \sum_{i = 1}^m x_i  = 1, \quad  \sum_{j = 1}^n y_j  = 1, \quad x,y \ge\boldsymbol{0}.
\end{align}
\end{subequations}
\end{minipage}
\hspace{3pt}

It is easy to see that any feasible mixed strategies $x,y$ will have objective function value less or equal to zero, as given player 2's (mixed) strategy, any pure strategy of player 1 can give payoff at most $p$, and given player 1's (mixed) strategy, any pure strategy of player 2 can give payoff at most $q$. This implies that any combination of pure strategies (i.e. any mixed strategy) of player 1 gives payoff at most $p$, and any combination of pure strategies of player 2 gives payoff at most $q$. Further, any Nash equilibrium $x^*,y^*$ has objective function value equal to zero, thus maximises the objective function. This is because players play best responses, and thus $p^* = {x^*}^\top Ay^*$ and $q^*={x^*}^\top By^*$. Importantly, \textit{only} the Nash equilibria have objective function value of zero. This is because for any optimal solution $(x^*,y^*,p^*,q^*)$ and any $(x,y)$ with $x \geq 0$, $\sum_{i=1}^m x_i =1$, $y \geq 0$, $\sum_{i=1}^n y_i =1$, $x^\top Ay^* \leq p^*$, ${x^*}^\top By \leq q^*$, and thus ${x^*}^\top Ay^* \leq p^*$, ${x^*}^\top By^* \leq q^*$. As a Nash equilibrium has objective function value of zero and is guaranteed to exist, the optimal value of the bilinear formulation must be zero (as it is non-positive). Thus ${x^*}^\top Ay^* = p^*$, ${x^*}^\top By^* = q^*$. This implies $x^\top Ay^* \leq {x^*}^\top Ay^*$, ${x^*}^\top By \leq {x^*}^\top By^*$.

In this work, we propose a multilinear feasibility program whose every feasible solution is a Nash equilibrium to a corresponding multi-player game with $n\geq2$ players. The formulation is based on an extension of the bilinear formulation to multi-player games. Although such an extension is straightforward, is has not been given much empirical analysis. We compare the multilinear feasibility formulation to established algorithms used to find an equilibrium in multi-player games. We find that our multilinear program is faster than a variety of algorithms in Gambit \cite{McKelvey2014}. 

\subsection{Multilinear formulation}
Let us define some notation. Let $n \geq 2$ be the number of players, and $[n] = \{1,\dots,n\}$ the set of players. Every player $i$ comes with a finite set of pure strategies $S_i$, with $|S_i| = n_i$. Let $\mathcal{S} = S_1 \times S_2 \times \dots \times S_n$ be the set of all $n$-tuples of pure strategy combinations of all players. We will further denote by $\mathcal{S}_{-i} = S_1 \times \dots \times S_{i-1} \times S_{i+1} \times \dots \times S_n$ the set of all pure strategy tuples of all players except $i$. Let $\bm{s} = (s_1,\dots,s_n) \in \mathcal{S}$ be a pure strategy tuple of all players and $\mathbf{\hat{s}} = (s_1,\dots,s_{i-1},s_{i+1},\dots,s_n) \in \mathcal{S}_{-i}$ be a pure strategy tuple of all players other than $i$. We define payoff matrix $A_i: \mathbf{s} \in \mathcal{S} \rightarrow \mathbb{R}$ for player $i$. As an example, if we had three players, $A_1[s_1,s_2,s_3]$ denotes player 1's payoff when player 1 plays pure strategy $s_1$, player 2 plays pure strategy $s_2$ and player 3 plays pure strategy $s_3$. Likewise, $A_2[s_1,s_2,s_3]$ and $A_3[s_1,s_2,s_3]$ denote player 2 and 3's payoff for the strategy combination $(s_1,s_2,s_3) \in \mathcal{S}$. For pure strategy $s$ of player $i$ and pure strategies $(s_1,\dots,s_{i-1},s_{i+1},\dots,s_n) = \mathbf{\hat{s}} \in \mathcal{S}_{-i}$ of the other players, we write $A_i[s,\mathbf{\hat{s}}]$ to denote the payoff of player $i$ when player $i$ plays pure strategy $s \in \mathcal{S}_i$ and the other players play pure strategies $\mathbf{\hat{s}} \in \mathcal{S}_{-i}$. For every player $i$, we define strategy vector $\mathbf{x}^i \in \mathbb{R}^{n_i}$, with $\mathbf{x}^i \geq \mathbf{0}$ and $\sum_{s \in S_i} x^i_s = 1$. $\mathbf{x}^i$ is a probability distribution over player $i$'s pure strategies, and thus a \textit{mixed strategy}. Note that any pure strategy is also a mixed strategy\footnote{When we refer to (mixed) strategies or a (mixed) Nash equilibrium, this includes pure strategies or a pure Nash equilibrium.}. Let $\mathbf{x} = (\mathbf{x}^1, \dots, \mathbf{x}^n)$ be a mixed strategy profile of all players, and $\mathbf{x}^{-i} = (\mathbf{x}^1, \dots, \mathbf{x}^{i-1}, \mathbf{x}^{i+1}, \mathbf{x}^n)$ be a mixed strategy profile of all players other than $i$. The product term $\prod_{s_j \in \mathbf{\hat{s}}} x_{s_j}^j$ for $\mathbf{\hat{s}} \in \mathcal{S}_{-i}$ denotes the combined probability of all players except $i$ to play the pure strategy tuple $\mathbf{\hat{s}} \in \mathcal{S}_{-i}$. As an example, if we have three players, player 2 has pure strategies $s_{2,1}$ and $s_{2,2}$ and player 3 has pure strategies $s_{3,1}$ and $s_{3,2}$, then $\mathcal{S}_{-1} = \{(s_{2,1},s_{3,1}),(s_{2,1},s_{3,2}),(s_{2,2},s_{3,1}),(s_{2,2},s_{3,2}) \}$. If player 2 plays $s_{2,1}$ with probability 1/2 and player 3 plays $s_{3,1}$ with probability 1/4, then $\prod_{s_j \in (s_{2,1},s_{3,1})} = 1/2*1/4$, $\prod_{s_j \in (s_{2,1},s_{3,2})} = 1/2*3/4$, $\prod_{s_j \in (s_{2,1},s_{3,1})} = 1/2*1/4$, $\prod_{s_j \in (s_{2,2},s_{3,2})} = 1/2*3/4$. Further, we define vector $\mathbf{p} \in \mathbb{R}^n$. $p^i$ corresponds to player i's highest expected payoff.
\begin{definition} \label{definition:equilibrium}
Let $\Gamma = \langle\{1,\dots,n\}, (S_i), (A_i) \rangle $ be a game with $n$, $S_i$, $A_i$, $\mathbf{x}^i$ defined as above. Let $\mathbf{x}^i \geq 0$ with $\sum_{s \in S_i} x^i_s = 1$ be a mixed strategy of player $i$. Then, $\mathbf{x^*} = (\mathbf{x^*}^1, \dots, \mathbf{x^*}^n)$ with $\mathbf{x^*}^i \geq 0$ and $\sum_{s \in S_i} {x^*}^i_s = 1$ for all players $i$ is a (mixed) Nash equilibrium if for all players $i$ and every mixed strategy $\mathbf{x}^i$, we have $\mathbb{E}\left[A_i[\mathbf{x^*}]\right] \geq \mathbb{E}\left[A_i[\mathbf{x}^i,\mathbf{x^*}^{-i}]\right]$.
\end{definition}

We now present the multilinear optimization formulation.

\LPblocktag{MLP1}{\label{fig:mlp1}}%
\begin{minipage}{12cm}
\begin{subequations}
\begin{alignat}{3}
	& \max_{\mathbf{x},\mathbf{p}} & \quad \sum_{i=1}^n \left( \sum_{\substack{(s,\mathbf{\hat{s}}) \\ \in S_i \times S_{-i}}} A_i[s,\mathbf{\hat{s}}] x_{s}^i \prod_{s_j \in \mathbf{\hat{s}}} x_{s_j}^j \right) & - \sum_{i=1}^n p^i & \label{eq:mlp1}\\
	& \text{   s.t.} & \quad \sum_{\mathbf{\hat{s}} \in S_{-i}} A_i[s,\mathbf{\hat{s}}] \prod_{s_j \in \mathbf{\hat{s}}} x_{s_j}^j & \leq p^i & \forall i \in [n], s \in S_i \label{eq:mlp2}\\
	& & \sum_{s \in S_i} x_s^i & = 1 & \forall i \in [n] \label{eq:mlp3}\\
	& & 0 \leq x_s^i & \leq 1 & \forall i \in [n], s \in S_i \label{eq:mlp4}
\end{alignat}
\end{subequations}
\end{minipage}

\begin{theorem} \label{multilinform}
A (mixed) strategy $(\mathbf{x}^1, \dots, \mathbf{x}^n)$ is a (mixed) Nash equilibrium of the n-player game $(A_1, \dots, A_n)$ if and only if there exist numbers $p^1,\dots,p^n$ such that $(\mathbf{x}^1, \dots, \mathbf{x}^n, p^1, \dots, p^n)$ is an optimal solution to the problem in \ref{fig:mlp1}.
\end{theorem}

It is easy to see that for two players, \ref{fig:mlp1} equals the bilinear formulation in \ref{fig:bimatrixform}, with $\mathbf{x}^1, \mathbf{x}^2$ instead of $\mathbf{x}, \mathbf{y}$, $p^1, p^2$ instead of $p,q$, and $A_1,A_2$ instead of $A,B$. Computational experiments on small instances reveal that the solver takes significant time to solve \ref{fig:mlp1} to optimality. However, further inspections reveal that it is more the verification of an optimal solution, rather than finding an optimal solution, that is the reason for this. Particularly, the solver finds a solution with objective function value 0 (i.e. a Nash equilibrium) relatively quickly, but spends a lot of time verifying that there is no feasible solution with objective function value larger than zero. However, as there cannot be a feasible solution with strictly positive objective value, it is sufficient for the solver to find a feasible solution with objective function value zero, instead of verifying that the upper bound to the optimisation program is zero. Thus, we reformulate the program into a feasibility program, for which the aim is to find a feasible solution for which the objective function (\ref{eq:mlp1}) of program \ref{fig:mlp1} is non-negative. As a strictly positive solution is not possible, any feasible solution to \ref{fig:mlp2} will have a value of zero, and thus be a Nash equilibrium. 

\begin{corollary}
Every feasible solution of \ref{fig:mlp2} is a Nash equilibrium.
\end{corollary}

\LPblocktag{MLP2}{\label{fig:mlp2}}%
\begin{minipage}{12cm}
\begin{subequations}
\begin{alignat}{3}
	& \max_{\mathbf{x},\mathbf{p}} & \quad 0 \\
	& \text{   s.t.} & \sum_{\mathbf{\hat{s}} \in S_{-i}} A_i[s,\mathbf{\hat{s}}] \prod_{s_j \in \mathbf{\hat{s}}} x_{s_j}^j & \leq p^i & \quad \forall i \in [n], s \in S_i \\
	& & \omit\rlap{$\displaystyle \sum_{i=1}^n \left( \sum_{\substack{(s,\mathbf{\hat{s}}) \\ \in S_i \times S_{-i}}} A_i[s,\mathbf{\hat{s}}] x_{s}^i \prod_{s_j \in \mathbf{\hat{s}}} x_{s_j}^j \right) - \sum_{i=1}^n p^i \geq 0$} \\
	& & \sum_{s \in S_i} x_s^i & = 1 & \forall i \in [n] \\
	& & 0 \leq x_s^i & \leq 1 & \forall i \in [n], s \in S_i 
\end{alignat}
\end{subequations}
\end{minipage}

\subsection{Mixed-integer formulations}
For a two-player game, four mixed-integer formulations whose solutions are equivalent to a Nash equilibrium in a two-player game were given in \cite{Sandholm2005}. We generalize these formulations to multi-player games. The notation we use is similar to the notation introduced in the multilinear formulations. Further, we introduce $U^i = \max_{s^l,s^h \in S_i, \mathbf{\hat{s}^l},\mathbf{\hat{s}^h}\in S_{-i}} A_i[s^h,\mathbf{\hat{s}^h}] - A_i[s^l,\mathbf{\hat{s}^l}]$ be the maximum difference of any two payoffs of player $i$ for any pure strategies of all players. 

We have four mixed-integer multilinear formulations, of which one is a feasibility program and three are optimisation programs. All programs have five sets of variables. $x_{s}^i, r_{s}^i, u_{s}^i, \overline{u}^i$ are real variables, and $b_{s}^i$ is binary. The MIMLPs have the same interpretation and range of values for variables $x_{s}^i \geq 0$, $\overline{u}^i, u_{s}^i, r_{s}^i \geq 0$, further they also come with constraints (\ref{eq:mip_2}), (\ref{eq:mip_3}), (\ref{eq:mip_4}), (\ref{eq:mip_5}) (which are mostly such that variables $x_{s}^i, \overline{u}^i, u_{s}^i, r_{s}^i$ are defined as desired). $x_{s}^i$, for all players $i \in [n]$ and all pure strategies $s \in S_i$ of player $i$, denotes the probability with which player $i$ plays pure strategy $s$. Hence, variables $x_{s}^i$ give us the mixed strategy played by each player. In order to be valid strategies, all pure strategies of a player must be played with non-negative probability (Eq. \ref{eq:mip_8}) and sum up to one (Eq. \ref{eq:mip_2}), for all players. $\overline{u}^i$ denotes the highest utility player $i$ can achieve by playing any strategy, given the other players mixed strategies. $u_{s}^i$ is the expected utility of player $i$ of playing pure strategy $s$, given the other players play their (mixed) strategies (Eq. \ref{eq:mip_3}). Naturally, $\overline{u}^i \geq u_{s}^i$ (Eq. \ref{eq:mip_4}). $r_{s}^i = \overline{u}^i - u_{s}^i$ (Eq. \ref{eq:mip_5}) is the regret of player $i$ of playing pure strategy $s$. It is defined as the difference of the highest utility of any strategy for the player to the utility of playing strategy $s$, given the other players' mixed strategies. By definition, the regret of any pure strategy must be non-negative (\ref{eq:mip_8}). Further, in any Nash equilibrium, every pure strategy that is played with strictly positive probability must have zero regret. If there was a pure strategy which the player plays and that has positive regret, the player can increase their payoff by putting more probability on a pure strategy with no regret and putting less probability on the pure strategy with regret. Hence, it would not be a Nash equilibrium. 

The meaning of binary variables $b_{s}^i$ is different in all formulations, with not all constraints of MIMLP 1 regarding this variable (Eq. \ref{eq:mip_6}, \ref{eq:mip_7}) present in MIMLP 2,3,4. In formulation 1, if $b_{s}^i$ is 1, strategy $s$ of player $i$ is not played, hence $x_{s}^i = 0$. If $b_{s}^i=0$, the probability on strategy $s$ is allowed to be positive, however the regret of the strategy must be zero. (\ref{eq:mip_6}) ensures that $b_{s}^i$ can only be set to 1 if zero probability is on $s$. Further, (\ref{eq:mip_7}) ensures that $b_{s}^i$ can only be set to zero if the strategy's regret is zero (if $b_{s}^i=1$, this constraint does not restrict any variable, as $r_{s}^i \leq U^i$ by definition).

\begin{proposition}	\label{mip1}
The set of feasible solutions to \ref{fig:mip1} is precisely the set of Nash equilibria for the corresponding multi-player game.
\end{proposition}

\LPblocktag{MIMLP1}{\label{fig:mip1}}%
\begin{minipage}{12cm}
\begin{subequations}
\begin{align}%{3}
 & \min & 0 & \label{eq:mip_1} \\
   & \text{   s.t.} & \sum_{s \in S_i} x_{s}^i & = 1 & \forall i \in [n] \label{eq:mip_2} \\
   & & u_{s}^i & = \sum_{\mathbf{\hat{s}} \in S_{-i}} \prod_{s_j \in \mathbf{\hat{s}}} x^j_{s_j} A_i[s, \mathbf{\hat{s}}] & \quad \forall i \in [n], \forall s \in S_i \label{eq:mip_3}\\
   & & \overline{u}^i & \geq u_{s}^i & \forall i \in [n], \forall s \in S_i \label{eq:mip_4}\\
   & & r_{s}^i & = \overline{u}^i - u_{s}^i & \quad \forall i \in [n], \forall s \in S_i \label{eq:mip_5}\\
   & & x_{s}^i & \leq 1 - b_{s}^i & \forall i \in [n], \forall s \in S_i \label{eq:mip_6}\\
   & & r_{s}^i & \leq U^ib_{s}^i & \forall i \in [n], \forall s \in S_i \label{eq:mip_7}\\
   & & x_{s}^i, r_{s}^i & \geq 0, \; u^i_s,\bar{u}_i \in \mathbb{R} & \forall i \in [n], \forall s \in S_i \label{eq:mip_8}\\
   & & b_{s}^i & \in \{0,1\} & \forall i \in [n], \forall s \in S_i \label{eq:mip_9}
\end{align}
\end{subequations}
\end{minipage}
\hspace{3pt}

\ref{fig:mip1} is a feasibility program, for which only Nash equilibria are feasible solutions. \ref{fig:mip2}, \ref{fig:mip3}, \ref{fig:mip4} have larger feasible regions, as pure strategies with positive probability are allowed to have positive regret, and pure strategies with positive regret are allowed to be played with positive probability. The formulations minimize a penalty, and it is only Nash equilibria for which the penalty is minimal. Thus, only Nash equilibria are optimal solutions. The advantage of these formulations is that, since finding a Nash equilibrium is assumed to be computationally intractable, these formulations can be used to stop the program before an equilibrium has been calculated, and thus give solutions which are close to an equilibrium, also called approximate equilibria. However, it is more difficult with these formulations to find a specific equilibrium among all equilibria, rather than just an arbitrary equilibrium. 

\ref{fig:mip2} penalises the regret of a pure strategy that is played with positive probability in the objective function, and thus for optimal solutions, the regret of pure strategies with positive probability is zero. \ref{fig:mip3} penalises the probability placed on pure strategies with positive regret, and thus optimal solutions will have zero probability on pure strategies with positive regret. \ref{fig:mip4} combines the normalised regret and the probability as a penalty, and the solver can choose whether the regret or the probability should be minimized. As \cite{Sandholm2005} noted, these formulations can be used to find approximate Nash equilibria.

\ref{fig:mip2} aims to minimize the regret of pure strategies that are played with positive probabilities. Particularly, the regret of a pure strategy played with positive probability serves as a penalty to the objective function. This is done by introducing variable $f_{s}^i$ for all $i \in [n], s \in S_i$, which represents a pure strategy's regret if the strategy has positive probability and zero otherwise. 

\begin{proposition} \label{mip2}
The set of Nash equilibria minimizes the objective function of \ref{fig:mip2}.
\end{proposition}

\LPblocktag{MIMLP2}{\label{fig:mip2}}%
\begin{minipage}{12cm}
\begin{subequations}
\begin{alignat}{3}
  & \min & & \sum_{i = 1}^n \sum_{s \in S_i} f_{s}^i - U^i b_{s}^i \label{eq:mip_10} \\
   & \text{   s.t.} & & (\ref{eq:mip_2}) - (\ref{eq:mip_6}),(\ref{eq:mip_8}),(\ref{eq:mip_9}) \nonumber \\
   & & f_{s}^i & \geq r^i_{s} & \forall i \in [n], \forall s \in S_i \label{eq:mip_11} \\
   & & f_{s}^i & \geq U^ib^i_{s} & \quad \forall i \in [n], \forall s \in S_i \label{eq:mip_12}
\end{alignat}
\end{subequations}
\end{minipage}
\hspace{3pt}

\ref{fig:mip3} is similar to MIMLP 2, however instead of minimising the regret of pure strategies played with positive probability, the probabilities of pure strategies with positive regret is minimized. To do so, variables $g_{s}^i$ are introduced, which are set such that a strategy's penalty in the objective function is zero if the strategy's regret is zero, and $x_{s}^i$ (the probability with which it is played) otherwise. The set of Nash equilibria minimizes the objective, as strategies with positive regret are not played. 

\begin{proposition} \label{mip3}
The set of Nash equilibria minimizes the objective function of \ref{fig:mip3}.
\end{proposition}

\LPblocktag{MIMLP3}{\label{fig:mip3}}%
\begin{minipage}{12cm}
\begin{subequations}
\begin{alignat}{3}
  & \min & & \sum_{i = 1}^n \sum_{s \in S_i} g_{s}^i - (1-b_{s}^i) \label{eq:mip_13} \\
   & \text{   s.t.} & & (\ref{eq:mip_2}) - (\ref{eq:mip_5}),(\ref{eq:mip_7}) - (\ref{eq:mip_9}) \nonumber \\
   & & g_{s}^i & \geq x_{s}^i & \forall i \in [n], \forall s \in S_i \label{eq:mip_14} \\
   & & g_{s}^i & \geq 1-b_{s}^i & \quad \forall i \in [n], \forall s \in S_i \label{eq:mip_15}
\end{alignat}
\end{subequations}
\end{minipage}
\hspace{3pt}

\ref{fig:mip4} combines MIMLP 2 and MIMLP 3. Instead of penalising all pure strategies' regret (MIMLP 2) or penalising all pure strategies' probabilities if they have positive regret (MIMLP 3), this formulation lets the solver decide for each pure strategy whether to penalise the regret or the probability. The penalised regret is expressed with variables $f_{s}^i$, the penalised probabilities are expressed with variables $g_{s}^i$. When using both the regret and the probabilities, the regret must be normalised, as the probability of a pure strategy is between zero and one, but a pure strategy's regret can generally be larger than one. Hence, $f_{s}^i$ uses normalised regret $r_s^i/U^i$, which is between zero and one.

\begin{proposition} \label{mimlp4}
The set of Nash equilibria minimizes the objective function of \ref{fig:mip4}.
\end{proposition}

\LPblocktag{MIMLP4}{\label{fig:mip4}}%
\begin{minipage}{12cm}
\begin{subequations}
\begin{alignat}{3}
  & \min & & \sum_{i = 1}^n \sum_{s \in S_i} f_{s}^i + g_{s}^i \label{eq:mip_16}\\
   & \text{   s.t.} & & (\ref{eq:mip_2}) - (\ref{eq:mip_5}),(\ref{eq:mip_8}),(\ref{eq:mip_9}) \nonumber \\
   & & f_{s}^i & \geq r_{s}^i/U^i & \forall i \in [n], \forall s \in S_i \label{eq:mip_17} \\
   & & f_{s}^i & \geq b_{s}^i & \quad \forall i \in [n], \forall s \in S_i \label{eq:mip_18} \\
   & & g_{s}^i & \geq x_{s}^i & \quad \forall i \in [n], \forall s \in S_i \label{eq:mip_19} \\
   & & g_{s}^i & \geq 1-b_{s}^i & \quad \forall i \in [n], \forall s \in S_i \label{eq:mip_20}
\end{alignat}
\end{subequations}
\end{minipage}
\hspace{3pt}

\subsection{Continuous and feasibility formulations}
For potential performance improvements of the mixed-integer multilinear programs, we further give continuous as well as feasibility formulations for the MIMLPs. Particularly, for all MIMLPs, we introduce continuous formulations\footnote{We note that due to this, the formulations are no longer mixed-integer, however we will still refer to MIMLP(C), to make clear that they belong to the respective MIMLP} MIMLP1(C), MIMLP2(C), MIMLP3(C), MIMLP4(C), for which constraint \ref{eq:mip_9}, i.e. constraints ($b_s^i \in \{0,1\}$) is replaced by $b_s^i = (b_s^i)^2$ (which implies $0 \leq b_s^i \leq 1$ and $b_s^i = 0$ or $b_s^i=1$). Thus, the continuous formulations are equivalent to the MIMLPs. The continuous formulation MIMLP1(C) for MIMLP1 is given in \ref{fig:mipmultilinear}, likewise MIMLP2(C), MIMLP3(C), MIMLP4(C) are simply MIMLP2, MIMLP3, MIMLP4 but constraint \eqref{eq:mip_9} replaced by \eqref{eq:mip_22}. 

\LPblocktag{MIMLP1(C)}{\label{fig:mipmultilinear}}%
\begin{minipage}{12cm}
\begin{subequations}
\begin{alignat}{3}
  & \min & & (\ref{eq:mip_1}) \nonumber \\
   & \text{   s.t.} & & (\ref{eq:mip_2}) - (\ref{eq:mip_5}),(\ref{eq:mip_8}) \nonumber \\
%   & & 0 \leq b_{s}^i & \leq 1 & \forall i \in [n], \forall s \in S_i \label{eq:mip_21} \\
   & & b_{s}^i & = (b_{s}^i)^2 & \quad \forall i \in [n], \forall s \in S_i \label{eq:mip_22} 
\end{alignat}
\end{subequations}
\end{minipage}
\hspace{3pt}

For MIMLP 2,3,4 we also introduce equivalent feasibility formulations MIMLP2(F), MIMLP3(F), MIMLP4(F), by introducing a constraint which requires the objective function of the respective MIMLP to be equal to the optimal value of the MIMLP. Particularly, MIMLP2 and MIMLP3 have optimal objective function of zero, and thus we introduce constraints (\ref{eq:mip_10}) = 0, i.e. $\sum_{i = 1}^n \sum_{s \in S_i} f_{s}^i - U^i b_{s}^i = 0$ (MIMLP2) and (\ref{eq:mip_13}) = 0, i.e. $\sum_{i = 1}^n \sum_{s \in S_i} g_{s}^i - (1-b_{s}^i) = 0$ for MIMLP3. MIMLP4 has optimal value $\sum_{i=1}^n |S_i|$, and thus we introduce constraint (\ref{eq:mip_16}) $= \sum_{i=1}^n |S_i|$, i.e. $\sum_{i = 1}^n \sum_{s \in S_i} f_{s}^i + g_{s}^i = \sum_{i=1}^n |S_i|$. For all feasibility formulations, the objective function is changed to 0. MIMLP3(F) is given in \ref{fig:mip3feasible}. 

Further, we introduce MIMLP2(C,F), MIMLP3(C,F), MIMLP4(C,F), which combine the continuous and feasibility formulations of MIMLP2,3,4, and are thus continuous multilinear formulations\footnote{MIMLP2(C,F), MIMLP3(C,F), MIMLP4(C,F) are thus no longer mixed-integer}. MIMLP3(C,F) is given in  \ref{fig:mip3similar}.

\LPblocktag{MIMLP3(F)}{\label{fig:mip3feasible}}%
\begin{minipage}{12cm}
\begin{equation*}
\begin{array}{l@{}c}
\text{min } & 0 \\
\text{s.t. } & (\ref{eq:mip_2}) - (\ref{eq:mip_5}),(\ref{eq:mip_7})-(\ref{eq:mip_9}),(\ref{eq:mip_14}),(\ref{eq:mip_15}) \\
%\text{subject to: } & (\ref{eq:mip_2}),(\ref{eq:mip_3}),(\ref{eq:mip_4}),(\ref{eq:mip_5}),(\ref{eq:mip_7}),(\ref{eq:mip_8}),(\ref{eq:mip_9}),(\ref{eq:mip_14}),(\ref{eq:mip_15}) \\
 & (\ref{eq:mip_13}) = 0 \\
  %& \text{minimize: } & & 0 \\
   %& \text{subject to: } & & (\ref{eq:mip_2}),(\ref{eq:mip_3}),(\ref{eq:mip_4}),(\ref{eq:mip_5}),(\ref{eq:mip_7}),(\ref{eq:mip_8}),(\ref{eq:mip_9}),(\ref{eq:mip_14}),(\ref{eq:mip_15}) \nonumber \\
   %& \text{subject to: } & & (\ref{eq:mip_2}) - (\ref{eq:mip_5}),(\ref{eq:mip_7})-(\ref{eq:mip_9}),(\ref{eq:mip_14}),(\ref{eq:mip_15}) \nonumber \\
   %& & (\ref{eq:mip_13}) & = 0 & \label{eq:mip_24} 
%\end{alignat}
\end{array}
\end{equation*}
\end{minipage}

\LPblocktag{MIMLP3(C,F)}{\label{fig:mip3similar}}%
\begin{minipage}{12cm}
\begin{equation*}
\begin{array}{l@{}c}
\text{min } & 0 \\
\text{s.t. } & (\ref{eq:mip_2}) - (\ref{eq:mip_5}),(\ref{eq:mip_7})-(\ref{eq:mip_8}),(\ref{eq:mip_14}),(\ref{eq:mip_15}) , (\ref{eq:mip_22}) \\
 & (\ref{eq:mip_13}) = 0 \\
%\end{alignat}
\end{array}
\end{equation*}
\end{minipage}

\begin{table}[htb]
\begin{threeparttable}
\centering
\begin{adjustbox}{width=\textwidth}
\begin{tabular}{l|cccccccccccccc}\toprule
MIMLP & 1 & 2 & 3 & 4 & 1(C) & 2(C) & 3(C) & 4(C) & 2(F) & 3(F) & 4(F) & 2(C,F) & 3(C,F) &4(C,F) \\ \midrule
Feasibility program & \checkmark &  &  &  & \checkmark &  &  &  & \checkmark & \checkmark & \checkmark & \checkmark & \checkmark & \checkmark \\
Optimality program &  & \checkmark & \checkmark & \checkmark &  & \checkmark & \checkmark & \checkmark &  &  &  &  &  &  \\
Continuous &  &  &  &  & \checkmark & \checkmark & \checkmark & \checkmark &  &  &  & \checkmark & \checkmark & \checkmark \\
Mixed-integer & \checkmark & \checkmark & \checkmark & \checkmark &  &  &  &  & \checkmark & \checkmark & \checkmark &  &  &  \\ \bottomrule
\end{tabular}
\end{adjustbox}
\end{threeparttable}
\caption{Overview of all mixed-integer multilinear formulations}
\end{table}

\section{Computational Experiments} \label{sec:computations}
All experiments are run on a MacBook Pro with 8GB RAM and Intel i5 CPU. Multilinear and mixed-integer formulations are implemented in AMPL \cite{Ampl2002}. We use BARON 21.1.13 \cite{BARON} as the solver, which uses FilterSD and FilterSQP as non-linear subsolvers. As the multilinear formulation MLP2 in \Cref{fig:mlp2} is much faster than any of the MIMLPs (see \Cref{fig:mimlp-graph}), we decide to only compare MLP2 against common algorithms for multi-player games. The MIMLPs do not seem to give better performance than existing algorithms, and hence the analysis of those is focused on comparing the MIMLPs to each other, to get an understanding which MIMLP formulation is best. Thus, the experiments consist of two parts: 
\begin{enumerate}
\item a comparison of MLP2 with common algorithms in \texttt{Gambit} \cite{McKelvey2014} (results in \Cref{fig:multilinear-algorithms-numerical}),
\item a comparison of the different MIMLPs (results in \Cref{fig:mimlp-graph}).
\end{enumerate}
All games are instanced in GAMUT \cite{Gamut2004} and have integer payoffs. We focus on \textit{random games} and \textit{covariance games} with negative covariance, as previous work \cite{Sandholm2005,Porter2008} indicates that covariance games with negative covariance are challenging to solve experimentally for a variety of algorithms as they tend to only have few equilibria with small support size. We refer to a covariance game with $n$ players and $|S_i|$ actions per player and covariance $\rho$ as CG(n,$|S_i|$,$\rho$), and to a random game with $n$ players and $|S_i|$ actions per player as RG(n,$|S_i|$). For all games, we take the average of 10 randomly generated instances of that game, and if a method did not find a Nash equilibrium before the timeout (which, depending on the game, is 300 or 900 seconds), we add the timeout to the average. 

\Cref{fig:multilinear-algorithms-numerical} compares MLP2 to the simplicial subdivision method (SD), the global newton method (GN), and the quantal response equilibrium (QRE) in Gambit. The results can be summarised as follows: The simplicial subdivision algorithm is the slowest, and already small instances are sufficient for the algorithm to not find a Nash equilibrium in less than 15 minutes. The global Newton method, although fast on the instances for which it finds an equilibrium, in many instances terminates without giving an equilibrium back. This issues has been reported in different scenarios, see \cite{TurocyAlgorithm}, and in these cases, we put the timeout towards the average. The logit algorithm and the multilinear formulation have similar runtime for smaller instances, but for larger games, our formulations seems to be faster. Thus, to conclude, our algorithm is faster than the algorithms in Gambit we test it with, and can be an alternative.

\Cref{fig:mimlp-graph} presents the results for the MIMLPs and the reformulations. It should be pointed out than any of the MIMLPs takes much longer to find an equilibrium than MLP2, and thus none of the MIMLPs is suited to find an equilibrium for \textit{large} multi-player games. This is different to the mixed-integer formulations for two-player games, for which \cite{Sandholm2005} showed better performance on some instances than existing algorithms. Therefore, the analysis of the MIMLPs aims more to get an understanding what type of formulation is best to find a Nash equilibrium in a multi-player game, than to compare the MIMLPs to common algorithms.

First, the continuous formulations MIMLP2(C), MIMLP3(C), MIMLP4(C) of MIMLP2, MIMLP3, MIMLP4 don't give much performance improvement compared to MIMLP2,3,4. For MIMLP2 and MIMLP3, both the feasibility formulations MIMLP2(F) and MIMLP3(F) and the combined continuous and feasibility formulations MIMLP2(C,F) and MIMLP3(C,F) give better performance than MIMLP2 and MIMLP3, but whether MIMLP2(C,F) and MIMLP3(C,F) are better than MIMLP2(F) and MIMLP3(C,F) depends very much on the game. For MIMLP4, whether MIMLP4(F) or MIMLP4(C,F) are better than MIMLP4 depends on the game. Further, compared over all games, MIMLP1(C), i.e. the continuous formulation of the feasibility formulation MIMLP1 seems to give the best performance.  
\vspace{-5pt}
\section{Future Work} \label{sec:conclusion}
Further questions include using different nonlinear solvers for the multilinear formulation. The solver we use finds a Nash equilibrium faster than any of the other algorithms we compare it to, other solvers should only improve the performance of the multilinear feasibility formulation. We also propose generating hard-to-solve instances. Even though GAMUT \cite{Gamut2004} offers many different types of games, many of these are easy to solve even for large multi-player games. Covariant games are among the few types of games that are (relatively) difficult to solve in the game generator GAMUT, and therefore we particularly use these instances. However, due to this, there is not much variety in the hard-to-solve instances we can use. Recent work has focused on hard-to-solve instances for polymatrix games (see \cite{gamesapproximation2016} and \url{http://polymatrix-games.github.io}), and so more hard-to-solve instances is a direction to explore. \looseness=-1

\begin{table}
\begin{threeparttable}
\centering
\begin{tabular}{p{3cm}ccccp{0.1cm}r}\toprule
Instance & MLP2 & GN & SD & QRE & & \\ \midrule
CG(5,5,$\rho=-0.2$) & 2.35 & 810.09 & 900 & 1.9 & & average (in seconds) \\
& 100\% & 10\% & 0\% & 100\% & & percentage solved \\ 
& 2.53 & 0.91 & -- & 1.9 & & average on solved (in seconds) \\[15pt] %\hline
CG(3,10,$\rho=-0.2$) & 0.57 & 271.56 & 518.96 & 0.36 & & average (in seconds) \\
& 100\% & 70\% & 50\% & 100\% & & percentage solved \\ 
& 0.57 & 2.22 & 137.9 & 0.36 & & average on solved (in seconds) \\[15pt] %\hline
RG(5,5) & 2.23 & 540.57 & 632.98 & 2.08 & & average (in seconds) \\
& 100\% & 40\% & 40\% & 100\% & & percentage solved \\ 
& 2.23 & 1.43 & 232.45 & 2.08 & & average on solved (in seconds) \\[15pt] %\hline
RG(3,10) & 0.329 & 91.325 & 382.9 & 0.362 & & average (in seconds) \\
& 100\% & 90\% & 70\% & 100\% & & percentage solved \\
& 0.329 & 1.47 & 161.287 & 0.362 & & average on solved (in seconds) \\[15pt] %\hline
CG(5,10,$\rho=-0.2$) & 250.28 & 825.52 & 900 & 361.46 & & average (in seconds) \\
& 100\% & 10\% & 0\% & 100\% & & percentage solved \\
& 250.28 & 155.21 & -- & 361.46 & & average on solved (in seconds) \\[15pt] %\hline
RG(5,10) & 208.79 & 900 & 900 & 564.32 & & average (in seconds) \\
& 100\% & 0\% & 0\% & 90\% & & percentage solved \\
& 208.79 & -- & -- & 527.02 & & average on solved (in seconds) \\[15pt] %\hline
CG(5,10,$\rho=-0.1$) & 220.72 & 813.47 & 900 & 415.64 & & average (in seconds) \\
& 100\% & 10\% & 0\% & 90\% & & percentage solved \\
& 220.72 & 34.77 & -- &361.82  & & average on solved (in seconds) \\\bottomrule
\end{tabular}
\footnotesize The time is the average over 10 instances of this game in seconds - if no solution is found after the timeout of 15 minutes, the timeout is evaluated as time for the instance.\\
\end{threeparttable}
\caption{Comparison of multilinear feasibility program to state-of-the-art algorithms} \label{fig:multilinear-algorithms-numerical}
\end{table}

\begin{table}
\centering
\begin{threeparttable}
\begin{tabular}{p{2.5cm} p{1.5cm} p{1.5cm} p{2cm} p{2cm} p{1.5cm}}\toprule
Method & RG(3,5) & RG(3,10) & CG(3,5,-0.2) & CG(5,3,-0.2) & RG(5,3) \\ \cmidrule{2-6} 
& Time & Time & Time & Time & Time \\ \midrule
MIMLP1 & 8.41 & 229.86 & 107.4 & 122.73 & 112.5 \\
MIMLP1(C) & 2.876 & 231.57 & 41.4 & 47.96 & 66.3\\ 
MIMLP2 & 19.11 & 202.38 & 16.16 & 538.17 & 660.5 \\
MIMLP2(C) & 55.93 & 272.15 & 165.5 & 469.5 & 453.14 \\
MIMLP2(F) & 14.4 & 150.33 & 29.72 & 410.94 & 345.45 \\
MIMLP2(C,F) & 14.3 & 279.03 & 68.516 & 198.7 & 218.16 \\
MIMLP3 & 46.79 & 265.53 & 161.4 & 392.45 & 535.9 \\
MIMLP3(C) & 75.93 & 300 & 200.59 & 575.32 & 430.03 \\
MIMLP3(F) & 17.67 & 225.66 & 18.8 & 188.94 & 115.95 \\
MIMLP3(C,F) & 9.16 & 260.98 & 76.71 & 54.56 & 90.91 \\
MIMLP4 & 5.84 &  220.969 & 79.4 & 359.54 & 49.75  \\
MIMLP4(C) & 110.3 & 300 & 69.6 & 479.0 & 462.62 \\
MIMLP4(F) & 56.5 & 221.26 & 129.78 & 129.59 & 56.88 \\
MIMLP4(C,F) & 59.19 & 270.69 & 65.12 & 248.85 & 84.52 \\ \hline
MLP 2 & 0.03 & 0.36 & 0.035 & 0.12 & 0.09 \\ \bottomrule
\end{tabular}
\footnotesize The time is the average over 10 instances of this game in seconds - if no solution is found after the timeout, the timeout is evaluated as time for the instance \\
RG(3,5), RG(3,10), CG(3,5,-0.2): Timeout after 300 seconds [5 minutes] \\
CG(5,3,-0.2), RG(5,3): Timeout after 900 seconds [15 minutes] \\
\end{threeparttable}
\caption{MIMLP results} \label{fig:mimlp-graph}
\end{table}

\ifSEA
	\bibliography{../multilinearNash.bib}
	\appendix
\else
{
	\newrefcontext[sorting=nyt]
	\printbibliography
}
\begin{appendices}
\fi

\section{Correctness of the Proposed Formulations}

Here we present proofs for the claims made with regards to the formulations presented in this paper.

\begin{proof}[\textbf{Proof of \Cref{multilinform}}]
We first show that if $(\mathbf{\bar{x}}^1, \dots, \mathbf{\bar{x}}^n)$ is a Nash equilibrium to $(A_1,\dots,A_m)$, then there exist numbers $\bar{p}^1,\dots,\bar{p}^n$ such that $(\mathbf{\bar{x}}^1,\dots,\mathbf{\bar{x}}^n,\bar{p}^1,\dots,\bar{p}^n)$ is an optimal solution to the program in \ref{fig:mlp1}.
Assume that $(\mathbf{\bar{x}}^1, \dots, \mathbf{\bar{x}}^n)$ is a Nash equilibrium. For any feasible solution $(\mathbf{x}^1,\dots,\mathbf{x}^n,p^1,\dots,p^n)$ of \ref{fig:mlp1}, constraints (\ref{eq:mlp2}), (\ref{eq:mlp3}) imply 
\[ \sum_{s \in S_i} \left( x_s^i \sum_{\mathbf{\hat{s}} \in S_{-i}} A_i[s,\mathbf{\hat{s}}] \prod_{s_j \in \mathbf{\hat{s}}} x_{s_j}^j \right) \leq p^i\] for all $i \in [n]$.  
This implies (\ref{eq:mlp1}) $\leq 0$ for any feasible solution. Set 
\[ \bar{p}^i = \sum_{s \in S_i} \left(\bar{x}_s^i \sum_{\mathbf{\hat{s}} \in S_{-i}} A_i[s,\mathbf{\hat{s}}]\prod_{s_j \in \mathbf{\hat{s}}} \bar{x}_{s_j}^j \right)\] 
for every $i \in [n]$. We show that $(\mathbf{\bar{x}}^1, \dots, \mathbf{\bar{x}}^n, \bar{p}^1, \dots, \bar{p}^n)$ is feasible and optimal to \ref{fig:mlp1}. 

As $(\mathbf{\bar{x}}^1, \dots, \mathbf{\bar{x}}^n)$ is a Nash equilibrium, we have 
\[ \sum_{s \in S_i} \left( \bar{x}_s^i \sum_{\mathbf{\hat{s}} \in S_{-i}} A_i[s,\mathbf{\hat{s}}] \prod_{s_j \in \mathbf{\hat{s}}} \bar{x}_{s_j}^j \right) \geq \sum_{s \in S_i} \left( x_s^i \sum_{\mathbf{\hat{s}} \in S_{-i}} A_i[s,\mathbf{\hat{s}}] \prod_{s_j \in \mathbf{\hat{s}}} \bar{x}_{s_j}^j \right)\]
for all (mixed) strategies $\mathbf{x}^i \geq 0$ with $\sum_{s \in S_i} x_s^i = 1$. Choosing $\mathbf{x}^i = \mathbf{e}_k$, with $k \in \{1,\dots,|S_i|\}$, hence the unit vector with all zeros except one in the k-th component, we have
\[ \bar{p}^i = \sum_{s \in S_i} \left( \bar{x}_s^i \sum_{\mathbf{\hat{s}} \in S_{-i}} A_i[s,\mathbf{\hat{s}}] \prod_{s_j \in \mathbf{\hat{s}}} \bar{x}_{s_j}^j \right) \geq \sum_{\mathbf{\hat{s}} \in \mathcal{S}_{-i}} A_i[k,\mathbf{\hat{s}}] \prod_{s_j \in \mathbf{\hat{s}}} \bar{x}_{s_j}^j \quad \forall k \in \{1,\dots,|S_i|\},\] satisfying constraint \cref{eq:mlp2}. As we can apply this for all players $i \in [n]$ (and constraints \cref{eq:mlp3,eq:mlp4} hold as $(\mathbf{\bar{x}}^1,\dots,\mathbf{\bar{x}}^n)$ is a Nash equilibrium), it follows that $(\mathbf{\bar{x}}^{1}, \dots, \mathbf{\bar{x}}^{m}, \bar{p}^1, \dots, \bar{p}^n)$ is feasible.
Further, the objective function value is zero at the point $(\mathbf{\bar{x}}^1, \dots, \mathbf{\bar{x}}^n, \bar{p}^1, \dots, \bar{p}^n)$. As the objective function value is at most zero and $(\mathbf{\bar{x}}^1, \dots, \mathbf{\bar{x}}^n, \bar{p}^1, \dots, \bar{p}^n)$ is feasible, it follows that it is optimal.
 
To show that if $(\mathbf{\bar{x}}^1,\dots,\mathbf{\bar{x}}^n,\bar{p}^1,\dots,\bar{p}^n)$ is an optimal solution to \ref{fig:mlp1}, $(\mathbf{\bar{x}}^1,\dots,\mathbf{\bar{x}}^n)$ is a Nash equilibrium, we assume that $(\mathbf{\bar{x}}^1,\dots,\mathbf{\bar{x}}^n,\bar{p}^1,\dots,\bar{p}^n)$ indeed is optimal to \ref{fig:mlp1}. Since a Nash equilibrium exists in this game and has objective value of zero, and all feasible solutions have non-positive value, it follows that the objective value of $(\mathbf{\bar{x}}^1,\dots,\mathbf{\bar{x}}^n,\bar{p}^1,\dots,\bar{p}^n)$ must be zero. 
For any $\mathbf{x}^i \geq 0$ with $\sum_{s \in S_i} x_s^i =1$, for all players $i \in [n]$, constraints \cref{eq:mlp2}, \cref{eq:mlp3} imply 
\[ \sum_{s \in S_i} \left( x_s^i \sum_{\mathbf{\hat{s}} \in S_{-i}} A_i[s,\mathbf{\hat{s}}] \prod_{s_j \in \mathbf{\hat{s}}} \bar{x}_{s_j}^j \right) \leq \bar{p}^i \quad \forall i \in [n] .\] Particularly, \[ \sum_{s \in S_i} \left( \bar{x}_s^i \sum_{\mathbf{\hat{s}} \in S_{-i}} A_i[s,\mathbf{\hat{s}}] \prod_{s_j \in \mathbf{\hat{s}}} \bar{x}_{s_j}^j \right) \leq \bar{p}^i \quad \forall i \in [n] .\] As further objective value of $(\mathbf{\bar{x}}^1,\dots,\mathbf{\bar{x}}^n,\bar{p}^1,\dots,\bar{p}^n)$ is zero, we have 
\[ \sum_{s \in S_i} \left( \bar{x}_s^i \sum_{\mathbf{\hat{s}} \in S_{-i}} A_i[s,\mathbf{\hat{s}}] \prod_{s_j \in \mathbf{\hat{s}}} \bar{x}_{s_j}^j \right) = \bar{p}^i \quad \forall i \in [n] .\] Hence, 
\[ \sum_{s \in S_i} \left( x_s^i \sum_{\mathbf{\hat{s}} \in S_{-i}} A_i[s,\mathbf{\hat{s}}] \prod_{s_j \in \mathbf{\hat{s}}} \bar{x}_{s_j}^j \right) \leq \sum_{s \in S_i} \left( \bar{x}_s^i \sum_{\mathbf{\hat{s}} \in S_{-i}} A_i[s,\mathbf{\hat{s}}] \prod_{s_j \in \mathbf{\hat{s}}} \bar{x}_{s_j}^j \right) = \bar{p}^i\]
$\forall i \in [n], \forall \mathbf{x}^i \geq 0: \sum_{s \in S_i} x_s^i=1$. Therefore, $(\mathbf{\bar{x}}^1,\dots,\mathbf{\bar{x}}^n)$ is a Nash equilibrium, as the constraint states that given the other players mixed strategies $\mathbf{\bar{x}}^j$, no strategy of player $i$ can give higher payoff than strategy $\mathbf{\bar{x}}^i$, for all players.
\end{proof}

\begin{proof}[\textbf{Proof of \Cref{mip1}}]
For any player $i \in [n]$ and any pure strategy $s \in S_i$ of player $i$, $x_{s}^i$ denotes the probability with which player $i$ plays pure strategy $s$. Constraint \ref{eq:mip_2},\ref{eq:mip_8} guarantee $\bm{x}^i$ to be a valid mixed strategy for each player $i$, as all pure strategies are played with non-negative probability and sum up to one. Constraint \ref{eq:mip_3} defines the expected payoff $u_s^i$ of player $i$ of playing pure strategy $s$ (given the other players' mixed strategies), and \ref{eq:mip_4} defines the highest possible expected payoff $\overline{u}^i$ of any (mixed) strategy of player $i$ given the other players' (mixed) strategies. Constraint \ref{eq:mip_5},\ref{eq:mip_8} define the regret $r_s^i$ of player $i$ of playing pure strategy $s \in S_i$. The regret of a pure strategy is the difference of player $i$'s highest possible expected payoff $\overline{u}^i$ and $i$'s payoff of playing pure strategy $s$ and is non-negative. Constraint \ref{eq:mip_9} introduces binary variable $b_s^i$ for any pure strategy $s$ of any player $i$. Constraint \ref{eq:mip_6} requires that $b_{s}^i$ can only be set to one if player $i$ puts zero probability on pure strategy $s$. Further, constraint \ref{eq:mip_7} ensures that $b_{s}^i$ can only be set to zero if the strategy's regret is zero (if $b_{s}^i=1$, this constraint does not restrict any variable, as $r_{s}^i \leq U^i$ by definition). Thus, if $b_{s}^i$ is 1, strategy $s$ of player $i$ is not played, hence $x_{s}^i = 0$. If $b_{s}^i=0$, the probability on strategy $s$ is allowed to be positive, however the regret of the strategy must be zero. Hence, only pure strategies with zero regret can be played with positive probability, which is precisely the definition of a Nash equilibrium.
\end{proof}

\begin{proof}[\textbf{Proof of \Cref{mip2}}]
Constraints \ref{eq:mip_2},\ref{eq:mip_3},\ref{eq:mip_4},\ref{eq:mip_5},\ref{eq:mip_8},\ref{eq:mip_9} guarantee that $\bm{x}^i,\bm{r}^i,\bm{u}^i,\overline{u}^i$ are correctly defined for all players. Due to constraint (\ref{eq:mip_6}), $b_{s}^i$ can only be set to one if the probability on the pure strategy is zero. Then, due to minimising $f_{s}^i$ in (\ref{eq:mip_10}) and Equation (\ref{eq:mip_12}), $f_{s}^i$ must be set to $U^i$. In that case, $f_{s}^i$ and $U^i b_{s}^i$ cancel out in the objective, and hence strategies with zero probability have no penalty. If $b_{s}^i = 0$, $f_{s}^i$ equals $r_{s}^i$, due to minimising $f_{s}^i$ and Equation (\ref{eq:mip_11}), and as $U^i b_{s}^i=0$, the penalty of the pure strategy equals the regret of the strategy, and pure strategies that have no regret do not have a penalty. Thus, due to the objective function, it is encouraged to play pure strategies which have no regret, and to not play strategies with regret. Thus, any pure strategy will only contribute to the objective function if it has positive regret \textit{and} probability. The Nash equilibria minimize the objective function, with optimal objective of zero. As any pure strategy in a Nash equilibrium will either have zero probability (hence no penalty) or zero regret (hence no penalty), the objective function will equal zero. Solutions which do not equal a Nash equilibrium have higher objective value, as for some strategies, $f_{s}^i > 0$ (as $r_{s}^i$ and $U^i$ are non-negative).
\end{proof}

\begin{proof}[\textbf{Proof of \Cref{mip3}}]
We recall that because of constraint (\ref{eq:mip_7}), $b_{s}^i$ can only be set to zero if the strategy's regret $r_{s}^i$ is zero. By constraint (\ref{eq:mip_15}) and minimising $g_{s}^i$, if $b_{s}^i = 0$, then $g_{s}^i = 1$. Thus, $g_{s}^i$ and $1-b_{s}^i$ cancel out in the objective function and the penalty of strategy $s$ is zero. If $b_{s}^i = 1$, due to constraint (\ref{eq:mip_14}) and minimising $g_{s}^i$, $g_{s}^i = x_{s}^i$, and $1-b_{s}^i = 0$. Hence, the penalty of strategy $s$ equals $x_{s}^i$. Therefore, the probability a pure strategy is played with only contributes to the objective function if the strategy has positive regret. Nash equilibria minimize the objective function, and come with optimal value of zero. Constraint (\ref{eq:mip_6}) of MIMLP 1 (namely, $x_{s}^i \leq 1-b_{s}^i$) is no longer in this formulation, and it is possible to set $b_{s}^i=1$ even if some probability is placed on $s$. However, in a Nash equilibrium, $b_{s}^i$ will only be set to 1 if the probability on $s$ is indeed zero, as pure strategies with positive regret are not played.
\end{proof}

\begin{proof}[\textbf{Proof of \Cref{mimlp4}}]
Constraint (\ref{eq:mip_17}) demands that if $b_{s}^i = 0$, then $f_{s}^i = r_{s}^i / U^i$, which is at most 1. Further, due to (\ref{eq:mip_20}), $g_{s}^i=1$. If $b_{s}^i = 1$, then $f_{s}^i = 1$ (constraint (\ref{eq:mip_18})) and $g_{s}^i = x_{s}^i$ (constraint (\ref{eq:mip_19})), which is at most 1. Hence, $f_{s}^i + g_{s}^i$ is at least 1 for every pure strategy $s$, and additional penalties (either the normalised regret or the probability of the strategy) contribute to the objective function if a strategy has positive probability and positive regret. Any feasible solution that is not a Nash equilibrium has $f_{s}^i + g_{s}^i > 1$ for some strategies, as not all strategies have either no regret or zero probability. Nash equilibria minimize the objective function, with value of $\sum_{i=1}^n |S_i|$, as the normalised regret is zero, or the probability of strategy is zero. 
\end{proof}

\ifSEA
\else
\end{appendices}
\fi

\end{document}